\documentclass[12pt]{article}
\usepackage[dvips]{epsfig}
\usepackage{amsmath}
\usepackage{amsfonts}
\usepackage[psamsfonts]{amssymb}
 \usepackage{amsthm}
\usepackage{bm}
\usepackage{enumerate}
\usepackage[mathscr]{eucal}
\usepackage{inputenc}
\usepackage[T2A]{fontenc}
\theoremstyle{plain}
\newtheorem{theorem}{Theorem}[section]
\newtheorem{lemma}{Lemma}[section]
\newtheorem{definition}{Definition}[section]
\newtheorem{remark}{Remark}[section]
\newtheorem{question}{Question}[section]
\newtheorem{example}{Example}[section]

 \numberwithin{equation}{section}

\begin{document}

\title{\textbf{The truncated Fourier operator. I}
\footnotetext{\hspace*{-4.0ex}\textbf{Mathematics Subject
Classification: (2000).} Primary 47A38; Secondary 47B35, 47B06,
47A10.\endgraf \hspace*{-2.0ex}\textbf{Keywords:} Truncated
Fourier-Plancherel operator, normal operator, contractive
operator, Hilbert-Schmidt operator, trace class operator.} }
\author{Victor Katsnelson \and Ronny Machluf}
\date{\ }%

 \maketitle

 \abstract{Let \(\mathscr{F}\) be the one dimensional Fourier-Plancherel operator
and \(E\) be a subset of the real axis. The truncated Fourier
operator is the operator \(\mathscr{F}_E\) of the form
\(\mathscr{F}_E=P_E\mathscr{F}P_E\), where
\((P_Ex)(t)=\chi_E(t)x(t)\), and \(\chi_E(t)\) is the indicator
function of the set \(E\). In the presented first part of the
work, the basic properties of the operator \(\mathscr{F}_E\)
according to the set \(E\) are discussed.

 Among these properties
there are the following one. The operator \(\mathscr{F}_E\): 1.
has a  not-trivial null-space; 2. is strictly contractive; 3. is a
normal operator; 4. is a Hilbert-Schmidt operator; 5. is a trace
class operator.}%

\section{The truncated Fourier operator:\\ %
 definition, basic properties}
Let \(E\) be a measurable subset of the real axis \(\mathbb{R}\).
(The case \(E=\mathbb{R}\) is not excluded). For \(p\geq 1\), let
\(L^p(E)\) be the space of complex valued functions on \(E\)
satisfying the condition \(\int\limits_{E}|x(t)|^p\,dt<\infty\,.\)
We mainly deal with the case \(p=2\), but episodic the case
\(p=1\) is needed. The space \(L^2(E)\), provided by the standard
linear operations and the scalar product
\(\langle{}x\,,\,y\,\rangle_E\):
\begin{equation}%
\label{1}
\langle{}x,y{}\rangle_E=\int\limits_{E}\,x(t)\overline{y(t)}\,dt,\quad
\big(x,y\in{}L^2(E)\big)\,,
\end{equation}
 is a Hilbert space. The norm in \(L^2(E)\) is

 \begin{equation}
 \label{2}
 \Vert{}x\Vert_E=\sqrt{\langle{}x,x{}\rangle_E}\,.
\end{equation}

The Fourier operator \(\mathscr{F}\) is defined by the formula
\begin{equation}
\label{3} %
(\mathcal{F}x)(t)=
\frac{1}{\sqrt{2\pi}}\int\limits_{\mathbb{R}}e^{it\xi}\,x(\xi)\,d\xi\,,\quad{}(t\in\mathbb{R}).
\end{equation}
One of central facts of the Fourier transform theory is the
Parseval equality:
\begin{equation}
\label{4}%
 \Vert(\mathscr{F}x)\Vert^2_{\mathbb{R}}=\Vert{}x\Vert^2_{\mathbb{R}}\,\quad (\forall
 x\in{}L^2(\mathbb{R})).
\end{equation}
This means that the Fourier operator \(\mathscr{F}\) is an
isometric operator in \(L^2(\mathbb{R})\). The next central  fact
of the Fourier transform theory is that the Fourier operator
\(\mathscr{F}\) maps the space \(L^2(\mathbb{R})\) \textit{onto}
the whole space \(L^2(\mathbb{R})\), that is the \(\mathscr{F}\)
is an \textit{unitary operator} in \(L^2(\mathbb{R})\). Moreover,
the inverse operator \(\mathscr{F}^{-1}\) is determined by the
formula
\begin{equation}%
\label{5}%
(\mathscr{F}^{-1}x)(t)=\frac{1}{\sqrt{2\pi}}\int\limits_{\mathbb{R}}e^{-it\xi}\,x(\xi)\,d\xi\,,\quad{}
(\,t\in\mathbb{R}).
\end{equation}%
\begin{remark} The integral in the right side hand of \eqref{3}
is a Lebesgue integral. It is well defined only if
\(x\in{}L^1(\mathbb{R})\). If
\(x\in{}L^1(\mathbb{R})\cap{}L^2(\mathbb{R})\), then both the
integral is well defined and the Parseval equality \eqref{4}
holds. Thus,  the operator \(\mathscr{F}\) can be
 defined originally by \eqref{3} only
 for
\(x\in{}L^1(\mathbb{R})\cap{}L^2(\mathbb{R})\). The set of such
\(x\) is dense in \(L^2(\mathbb{R})\). Since the operator
\(\mathscr{F}\) acts isometrically on this set, it can be extended
by the continuity on the whole space \(L^2(\mathbb{R})\). The same
is related to the operator which appears in \eqref{5}.
\end{remark}

\textsf{In this paper we deal with the truncated Fourier
operator.}

\begin{definition}
Let \(E\) be a measurable subset of the real axis,
\begin{math}
 0<m(E)\leq{}\infty\,.
\end{math}
 The operator
\(\mathscr{F}_E:\,L^2(E)\to{}L^2(E)\), is defined as
\begin{equation}
\label{6}%
 (\mathscr{F}_Ex)(t)=
\frac{1}{\sqrt{2\pi}}\int\limits_{E}e^{it\xi}\,x(\xi)\,d\xi\,,\quad{}(t\in{}E).
\end{equation}
The operator \(\mathscr{F}_Ex\) is said to be \textsf{the
truncated Fourier operator}, or in more detail, \textsf{the
 Fourier operator truncated on the set \textit{E}}.%
 \end{definition}

\begin{remark}
If the set \(E\) is a set of finite Lebesgue measure:
\(\int_Edt<\infty\), then \(L^2(E)\subset{}L^1(E)\) and the
integral in \eqref{6} is well defined for every \(x\in{}L^2(E)\).
\end{remark}

The operator \(\mathscr{F}_E^{\,\ast}\), which is the adjoint
operator to \(\mathscr{F}_E\) with respect to the scalar product
\eqref{1}, is
\begin{equation}
\label{8}%
(\mathcal{F}_E^{\,\ast}x)(t)=
\frac{1}{\sqrt{2\pi}}\int\limits_{E}e^{-it\xi}\,x(\xi)\,d\xi\,,\quad{}(t\in{}E).
\end{equation}

\begin{remark}
\label{OtDe}%
The operator \(\mathscr{F}_E\), acting in \(L^2(E)\), may be
naturally identified with the operator \(P_E\mathscr{F}P_E\),
acting in \(L^2(\mathbb{R})\), where
 \begin{equation}%
 \label{CaPr}%
 (P_Ex)(t)=\chi_E(t)x(t),
 \end{equation}%
\begin{equation}%
\label{CaF}%
 \chi_E(t)=
\begin{cases}
 1,& t\in{}E\,,\\
0,& t\notin{}E\,.
\end{cases}
\end{equation}%
\end{remark}

\begin{lemma}
For any \(E\), the operator \(\mathscr{F}_E\) is a contractive
operator in \(L^2(E)\):
\begin{equation}%
\label{7}%
 \Vert\mathscr{F}_Ex\Vert^2_E\leq\Vert{}x\Vert^2_E\quad(\forall{}x\in{}L^2(E)).
\end{equation}
\end{lemma}
\begin{proof} Indeed, if
\begin{equation}
\label{p}
y(t)=\frac{1}{\sqrt{2\pi}}\int_{E}e^{it\xi}\,x(\xi)\,d\xi\,,\
(t\in\mathbb{R}),
\end{equation}
 then by the Parseval equality,
\(\int_{\mathbb{R}}|y(t)|^2=\int_E|x(t)|^2\), hence
\(\int_E|y(t)|^2\leq\int_E|x(t)|^2\).
\end{proof}
Thus, the inequalities hold
\begin{equation}%
\label{Co}%
 0\leq\mathscr{F}_E^\ast\mathscr{F}_E\leq{}I_E \textup{
and } 0\leq\mathscr{F}\mathscr{F}^\ast\leq{}I_E\,.
\end{equation}%
Here and in what follows \(I_E\) is the identity operator.

\begin{theorem}\label{IP} {\ }\\
\hspace*{2.0ex}\textup{1}\,.\,\,If \(\textup{mes}\,E>0\), then
no-one of the operators \(\mathscr{F}_E=0\) and
\(\mathscr{F}_E^\ast=0\) equals zero: there exists
\(x\in{}L^2(E)\) for which both
\begin{equation}%
\label{CE1}%
\mathscr{F}_Ex\not=0\quad\text{and}\quad
\mathscr{F}_E^\ast{}x\not=0\,.
\end{equation}
\\[0.0ex]
\hspace*{2.0ex}\textup{2}\,.\,\,If \,%
\(\textup{mes}(\mathbb{R}\setminus{}E)>0\), then no-one of the
operators \(\mathscr{F}_E\) and \(\mathscr{F}_E^\ast\) is
isometric: there exists \(x\in{}L^2(E)\) for which both
\begin{equation}%
\label{CE2}%
\Vert\mathscr{F}_Ex\Vert<\Vert{}x\Vert\quad\text{and}\quad%
\Vert\mathscr{F}_E^\ast{}x\Vert<\Vert{}x\Vert\,.
\end{equation}
\end{theorem}
\begin{proof}
Let \(t_0\in{}E\) is such that %
\begin{equation*}\lim_{n\to\infty}\frac{m\big(E\cup[t_0-1/n\,,\,t_0+1/n]\big)}{m\big([t_0-1/n\,,\,t_0+1/n]\big)}=1\,.
\end{equation*}
(Almost every point \(t_0\in{}E\) possesses  this property). Let
us set
\begin{equation*}
x_n(\xi)=\frac{1}{m\big([t_0-1/n\,,\,t_0+1/n]\big)}%
\begin{cases}
1, & \xi\in{}E\cap[t_0-1/n\,,\,t_0+1/n]\\
0, & \xi\in{}E\setminus[t_0-1/n\,,\,t_0+1/n]
\end{cases}\,.
\end{equation*}
It is clear that \(x_n\in{}L^2(E)\,\,\forall{}n\), and that
\begin{equation*}
\int\limits_{E}x_n(\xi)e^{\pm{}it\xi}\,d\xi\to{}e^{\pm{}i\,t\,t_0}\
\ \textup{as}\,\,n\to\infty,\quad\textup{the limit is locally
uniform on}\ \mathbb{R}\,.
\end{equation*}
If \(n\) is large enough, then
\(\int\limits_{\mathbb{R}\setminus{}E}|\mathscr{F}x_n)(t)|^2\,dt>0,
\ \int\limits_{\mathbb{R}\setminus{}E}|(\mathscr{F}^\ast%
x_n)(t)|^2\,dt>0\). Because \(x_n\) vanishes outside of \(E\),
\(\int\limits_{E}|x_n(t)|^2\,dt=\int\limits_{\mathbb{R}}|x_n(t)|^2\,dt\).
By the Parseval equality,
\(\int\limits_{\mathbb{R}}|\mathscr{F}x_n)(t)|^2\,dt=
\int\limits_{\mathbb{R}}|(\mathscr{F}^\ast{}x_n)(t)|^2\,dt=
\int\limits_{E}|x_n(t)|^2\,dt\).
 Thus, the inequalities  \eqref{CE2} hold for \(x=x_n\) if \(n\) is large enough.
\end{proof}

It is clear that if the set \(E\) is bounded, then then the
inequalities \eqref{CE1}, \eqref{CE2} hold for \textit{any}
\(x\in{}L^2(E),\,x\not=0\)\,. Indeed, given
\(x\in{}L^2(E),\,x\not=0\), let \(y\) is determined from \(x\)
according to \eqref{p}. Since the set \(E\) is bounded, \(y(t)\)
is an entire function of \(t\). Therefore the function \(y(t)\)
may vanish only in isolated points. In particular,
\(\int\limits_{E}|y(t)|^2>0,\,\int\limits_{\mathbb{R}\setminus{}E}|y(t)|^2>0\).
The first inequality means that \(\Vert\mathscr{F}x\Vert>0\), the
second one --- that \(\Vert\mathscr{F}x\Vert<\Vert{}x\Vert\).
 In
\cite[Proposition 5]{AmBe} it was shown that if
\(\textup{mes}\,E<\infty\) (the set \(E\) may be unbounded), then
the inequalities \eqref{CE2} holds for arbitrary
\(x\in{}L^2(E),\,x\not=0.\)

Actually, the much more stronger statement takes place.
\begin{theorem}%
\label{NKer}%
 If  \(\textup{mes}\,E<\infty\), then the inequalities
\label{NUP}
\begin{equation}
\label{BFU}%
 \Vert\mathscr{F}_Ex\Vert^2\leq{}\big(1-A^{-1}
e^{-A(\textup{mes}\,E)^2}\big)\Vert{}x\Vert^2\,,\quad{}
\Vert\mathscr{F}^{\ast}_Ex\Vert^2
\leq\big(1-A^{-1}e^{-A(\textup{mes}\,E)^2}\big)\Vert{}x\Vert^2%
\end{equation}
hold for every \(x\in{}L^2(E)\)\,. Here \(A, \ 1\leq{}A<\infty\),
is an absolute constant: it does not depend neither on \(x\), nor
on \(E\).
\end{theorem}

\begin{proof} In fact,
Theorem \ref{NKer} is a special case of the \textit{Nazarov
uncertainty principle}. In \cite{Naz}, F.L. Nazarov prove the
remarkable inequality
\begin{equation}
\label{NRE}%
 \int\limits_{\mathbb{R}}\vert{}y(t)\vert^2\,dt\leq
A\,e^{A(\textup{mes}\,E)(\textup{mes}\,F)}
\big(\int\limits_{\mathbb{R}\setminus{}E}\vert{}y(t)\vert^2\,dt+%
\int\limits_{\mathbb{R}\setminus{}F}\vert{}x(\xi)\vert^2\,d\xi\big)\,
\end{equation}
where \(x\in{}L^2(\mathbb{R})\) is an arbitrary functions, \(y\)
is the Fourier transform of \(x\): \({}y(t)=\frac{1}{\sqrt{2\pi}}
\int\limits_{\mathbb{R}}\,e^{it\xi}\,x(\xi)\,d\xi\), \(E\) and
\(F\) are arbitrary measurable subsets of \(\mathbb{R}\). If
\(F=E\)  and \(x\) is an arbitrary function vanishing outside of
\(E\), then the inequality  \eqref{NRE} takes the form
\[ \int\limits_{\mathbb{R}}\vert{}y(t)\vert^2\,dt\leq
A\,e^{A(\textup{mes}\,E)^2}
\big(\int\limits_{\mathbb{R}}\vert{}y(t)\vert^2\,dt-
\int\limits_{\mathbb{E}}\vert{}y(t)\vert^2\,dt\big)\,.
\]
 Invoking the Parseval
identity, we rewrite this inequality in the form
\[\int\limits_{\mathbb{E}}\vert{}x(t)\vert^2\,dt\leq
A\,e^{A(\textup{mes}\,E)^2}
\big(\int\limits_{\mathbb{E}}\vert{}x(t)\vert^2\,dt-
\int\limits_{\mathbb{E}}\vert{}y(t)\vert^2\,dt\big)\]

The latter inequality coincides with the first of the inequalities
\eqref{BFU}.
\end{proof}
\begin{remark}
\begin{upshape}
As it is stated below, the equality for the Hilbert-Schmidt norm
\(\Vert\mathscr{F}_E\Vert_{\mathfrak{S}_2}\) holds:
\(\Vert\mathscr{F}_E\Vert_{\mathfrak{S}_2}=\textup{mes}\,E\).
Since the Hilbert-Schmidt norm majorizes  the operator norm, the
inequalities hold
\begin{equation}
\label{SmM}
 \Vert\mathscr{F}_Ex\Vert\leq{}(\textup{mes}\,E)
\Vert{}x\Vert,\quad
 \Vert\mathscr{F}^\ast_Ex\Vert\leq{}(\textup{mes}\,E)\,
\Vert{}x\Vert\,,
\end{equation}
where \(E\) is an arbitrary measurable set and \(x\)is an
arbitrary function from \(L^2(E)\).

Both inequalities \eqref{BFU} and \eqref{SmM} are true. However
\eqref{SmM} is more precise for small values of
\(\textup{mes}\,E\), and \eqref{BFU} -- for large ones.
\end{upshape}
\end{remark}

\begin{remark}
\begin{upshape} If the set \(E\) is not just of a set of finite
measure, but a finite interval, the the estimate \eqref{BFU} can
be refined for large values of \(\textup{mes}\,E\). For
\(E=[-l,l]\), the largest eigenvalues \(\lambda_0(l)\) of the
operator \(\mathscr{F}^\ast_E\mathscr{F}_E\) coincides with the
squares of norm \(\Vert\mathscr{F}_E\Vert^2\).  For \(E=[-l,l]\),
the operator \(\mathscr{F}^\ast_E\mathscr{F}_E\) is the integral
operator in \(L^2([-l,l])\) of the form  %

\begin{equation}%
\label{AsB}
\big(\mathscr{F}^\ast_E\mathscr{F}_E\,x\big)(t)=%
\frac{1}{\pi}\,\int\limits_{[-l.l]}
\frac{\sin{}l(t-\tau)}{t-\tau}x(\tau)\,d\tau\,.
\end{equation}

 The asymptotic
behavior as \(l\to\infty\) of the eigenvalue \(\lambda_0(l)\) of
the integral operator \eqref{AsB}
 was
found by W.H.J.\,Fuchs in \cite{Fu}:
\[1-\lambda_0(l)\approx{}4\sqrt{\pi}l^{1/2}e^{-2l}
\,\,(=2\sqrt{2\pi}(\textup{mes}\,E)^{1/2}e^{-\textup{mes}\,E}.
\]
Thus, for \(E=[-l,l]\), the estimate  holds which is strongest
then the estimate \eqref{BFU}:

\begin{equation}
\Vert\mathscr{F}_E\,x\Vert\leq(1-A(\varepsilon)%
e^{-(1+\varepsilon)\textup{mes\,E}})\Vert{}x\Vert\,,\quad
\Vert\mathscr{F}^\ast_E\,x\Vert\leq(1-A(\varepsilon)%
e^{-(1+\varepsilon)\textup{mes\,E}})\Vert{}x\Vert\,,
\end{equation}
for every \(x\in{}L^2(E)\)\,. Here \(\varepsilon>0\) is arbitrary,
 and \(A(\varepsilon)<\infty\) for any \(\varepsilon>0\). The value
\(A(\varepsilon)\) does not depend on \(l\) and \(x\).
\end{upshape}
\end{remark}

\vspace{2.0ex}

Theorem \ref{IP} claims that if the set \(E\) is bounded, then the
null-spaces of each of operators
\(\mathscr{F}_E^\ast\mathscr{F}_E\),
\(\mathscr{F}^\ast_E\mathscr{F}^\ast_E\) are trivial\,---\,they
consist of zero-vector only.

Theorem \ref{NKer} implies that if \(\textup{mes}\,E<\infty\),
then the null-spaces of each of operators
\(I_E-\mathscr{F}^\ast_E\mathscr{F}_E\),
\(I-\mathscr{F}_E\mathscr{F}_E^\ast\) are trivial.

The following example shows that if  \(\textup{mes}\,E=\infty\) ,
then each of these null-spaces can be not only non-trivial, but
even  an infinite dimensional one.
\begin{example}
\begin{upshape}
Let \(K\subset\mathbb{R}\) be the interval:
\begin{equation}%
\label{bi}%
K=[-a,a], \textup{\ where\ } 0<a<\sqrt{\frac{\pi}{2}}.%
\end{equation}
 The
set \(E\) is a "periodic" systems of intervals:
\begin{equation}
\label{dfa}
E=\bigcup\limits_{p\in\mathbb{Z}}\big(K+p\sqrt{2\pi}\big)\,.
\end{equation}
Let \(u(t)\not\equiv0\) be a (smooth) function on \(\mathbb{R}\)
such that
\begin{equation}%
\label{df0}%
\text{supp}\,u\subseteq{}K\,.
\end{equation}
 The
function \(u(t)\) is representable in the form
\begin{equation}
\label{df1}%
u(t)=\int\limits_{-\infty}^{\infty}e^{it\xi}v(\xi)\,d\xi\,,
\end{equation}
where \(v(\xi)\) is a fast decaying function.
 Let
\begin{equation}
\label{df2}%
y(t)=\sum\limits_{p\in\mathbb{Z}}c_pu(t+p\sqrt{2\pi})\,,
\end{equation}
where \(\lbrace{}c_p\rbrace_{p\in\mathbb{Z}}\) be a summable
sequence. From \eqref{dfa}\,-\eqref{df2} it follows that
\begin{equation}
\label{dfb}%
\textup{supp}\,y\subseteq{}E\,.
\end{equation}
Moreover
\begin{equation}
\label{df3}%
y(t)=\int\limits_{-\infty}^{\infty}e^{it\xi}v(\xi)\,\varphi(\xi)\,d\xi,
\end{equation}
where
\begin{equation}
\label{df4}%
\varphi(\xi)=\sum\limits_{p\in\mathbb{Z}}c_pe^{ip\xi{}\sqrt{2\pi}},\quad
-\infty<\xi<\infty\,.
\end{equation}
The function \(\varphi\) is a periodic one:
\begin{equation}
\label{df5}%
 \varphi(\xi+\sqrt{2\pi})\equiv\varphi(\xi),\quad
-\infty<\xi<\infty\,.
\end{equation}
Let us invert the order of reasoning. Starting from a function
\(u(t)\) supported on \(K\), \eqref{df0}, and
\(\sqrt{2\pi}\)-periodic function \(\varphi(\xi)\), \eqref{df5},
we \textit{define} the function \(y(t)\) by \eqref{df3}, where
\(v(\xi)\) is determined from \(u\) by \eqref{df1}. Then the
equality \eqref{df2} holds, where
\(\lbrace{}c_p\rbrace_{p\in\mathbb{Z}}\) is the sequence of the
Fourier coefficient by the originally given function
\(\varphi\):\,\eqref{df4}. Let a \(\sqrt{2\pi}\) periodic function
\(\varphi\not\equiv{}0\) satisfy the condition
\begin{equation}
\label{df6}%
\textup{supp}\,\varphi\cap[\scriptstyle{-\sqrt{\pi/2},-\sqrt{\pi/2}}\textstyle]\subseteq{}K,
\end{equation}
where \(K\) is the same as before. Then
\begin{equation}
\label{df7}%
\textup{supp}\,v(\xi)\varphi(\xi)\subseteq{}E\,.
\end{equation}
If moreover the function \(\varphi\) is smooth, then the sequence
\(\lbrace{}c_p\rbrace_{p\in\mathbb{Z}}\), \eqref{df4}, is
summable. Thus the function \(y(t)\) is representable in the form
\eqref{p}, where \(x(\xi)=\sqrt{2\pi}v(\xi)\varphi(\xi)\),
\(\text{supp}\,x\subseteq{}E,\,\text{supp}\,y\subseteq{}E\),
therefore
\[\int\limits_E|y(t)|^2=\int\limits_{\mathbb{R}}|y(t)|^2=\int\limits_E|x(\xi)|^2d\xi\,,\]
i.e.
\begin{equation}%
\label{Iso}%
 \Vert\mathscr{F}_Ex\Vert^2=\Vert{}x\Vert^2\,.
\end{equation}%
Because of the freedom in the choice of \(u(t)\) and
\(\varphi(\xi)\), \textsf{the set of
\({\boldsymbol{x\in{}L^2(E)}}\) satisfying the condition
\eqref{Iso} is an infinite dimensional subspace of
\(\boldsymbol{L^2(E)}\).}

Let \(x_1(\xi)=x(\xi)e^{-ih\xi}\), where \(h\in\mathbb{R}\), and
\[y_1(t)=\frac{1}{\sqrt{2\pi}}\int_{E}e^{it\xi}\,x_1(\xi)\,d\xi\,,\quad
t\in\mathbb{R},\] Then
\(y_1(t)=y(t-h),\,\textup{supp}\,y_1=h+\textup{supp}\,y\). If
\(a<\frac{1}{2}\sqrt{\frac{\pi}{2}}\),\,then \(h\) can be chosen
such that \((E+h)\cap{}E=\emptyset\,.\) In this case,
\(y_1(t)=0\,\,\forall{}t\in{}E\), thus
\begin{equation}%
\label{Nul}%
 \mathscr{F}_Ex_1=0\,.
\end{equation}
 As before, \textsf{the set of
\(\boldsymbol{x_1\in{}L^2(E)}\) satisfying the condition
\eqref{Nul} is an infinite dimensional subspace of
\(\boldsymbol{L^2(E)}\).}

In this example, both
\begin{equation}%
\label{im}%
\textup{mes}(E)=\infty,\quad\textup{mes}(\mathbb{R}\setminus{}E)=\infty\,.
\end{equation}%
\end{upshape}
{\ }\hfill \(\Box\)
\end{example}

\begin{remark}%
\label{Exam}%
In \textup{\cite[Proposition 6]{AmBe}} it was shown that if a set
\(E\) satisfies the condition
\(\textup{mes}\,(\mathbb{R}\setminus{}E)<\infty\), then the set of
\(x\in{}L^2(E)\) satisfying the equality \eqref{Iso} is an
infinite dimensional subspace of \(L^2(E)\).
\end{remark}%

We recall that the operator \(A\) acting in a Hilbert space is
said to be \textit{normal} if
\[A^{\ast}A=AA^\ast\,,\]
where \(A^\ast\) is the operator adjoint to the operator \(A\).\\

Here and further
 \begin{equation}
\label{MiE}
-E=\lbrace{t\in\mathbb{R}: -t\in{}E\rbrace}%
 \end{equation}.

\begin{lemma}%
\label{NCo}%
 The truncated Fourier operator \(\mathscr{F}_E\) is normal
if and only if the
 equality
\begin{equation}
\label{NCoE}%
 \int\limits_{E\setminus(-E)}|y(t)|^2dt=\int\limits_{(-E)\setminus{}E}|y(t)|^2dt\,,
\end{equation}
holds for every \(y(t)\) of the form
\(y(t)=\frac{1}{\sqrt{2\pi}}\int_{E}e^{it\xi}\,x(\xi)\,d\xi\),
\(t\in\mathbb{R}\), where \(x\) runs over the whole space
\(L^2(E)\).
\end{lemma}
\begin{proof}
The condition
\(\mathscr{F}^\ast_E\mathscr{F}_E=\mathscr{F}_E\mathscr{F}^\ast_E\)
is equivalent to the condition: \textit{the equality
\(\Vert\mathscr{F}_Ex\Vert^2=\Vert\mathscr{F}^\ast_Ex\Vert^2\)
holds for every \(x\in{}L^2(E)\).} If \(x\in{}L^2(E)\), then
\((\mathscr{F}_{E}x)(t)=y(t),\,t\in{}E\,\), and
\((\mathscr{F}^{\ast}_{E}x)(t)=y(-t),\,t\in{}E\,\). Thus, the
equality
\(\Vert\mathscr{F}_Ex\Vert^2=\Vert\mathscr{F}^\ast_Ex\Vert^2\)
takes the form
\(\int\limits_{E}|y(t)|^2dt=\int\limits_{E}|y(-t)|^2dt\). The
latter equality is equivalent to the equality \eqref{NCoE}.
\end{proof}

\begin{definition}
\label{SiSe}%
 The set \(E\) is said to be \textit{symmetric} if
\begin{equation}
\label{SiCo}%
 \textup{mes}\,\Delta(E,-E)=0\,,
\end{equation}
where \(\Delta(E,-E)\) is the symmetric difference of the sets
\(E\) and \(-E\):
\[\Delta(E,-E)=(E\setminus(-E))\cup((-E)\setminus{}E).\]
\end{definition}
 Since
\((E\setminus(-E))\cap((-E)\setminus{}E)=\emptyset\), and %
\(\textup{mes}\,(E\setminus(-E))=\textup{mes}\,((-E)\setminus{}E)\),
the condition \eqref{SiCo} can be expressed in asymmetric form: %
\[\textup{mes}\,\Delta(E,-E)=0.\]

\begin{theorem}
\label{SiN}%
 If the set \(E\) is symmetric, then the operator
\(\mathscr{F}_E\) is a normal operator.
\end{theorem}
\begin{proof}
The theorem is an evident consequence of Lemma \ref{NCo}: the
expressions in both sides of \eqref{NCoE} are equal because both
of them vanish.
\end{proof}%

\begin{question}
Let the operator \(\mathscr{F}_E\) be normal. Is the set \(E\)
symmetric?
\end{question}

We can not answer this question in full generality. However, under
some extra condition imposed on the set \(E\) the answer  to this
question is affirmative.

The set \(S,\ S\subset{}\mathbb{R}\), is said to be
\textit{bounded,\,\,bounded from below, and bounded from above}
respectively, if \(S\) is contained respectively in some bounded
interval \([a,b]\), bounded from above interval \([a,+\infty)\) or
bounded from below interval \((-\infty,b]\), where \(a,b\) are
some finite numbers. (In the first case, \(a<b\).) The set \(S,\
S\subset{}\mathbb{R}\), is said to be \textit{semi-bounded}, if
\(S\) is either bounded from above, or is bounded from below. (In
particular, every bounded set is semi-bounded).
\begin{theorem}
\label{NiS}%
Assume that the following two conditions are satisfied:\\
\hspace*{2.0ex}\,\,\textup{1.} The operator \(\mathscr{F}_E\) is
normal;\\
 \hspace*{2.0ex}\,\,\textup{2.} The set \(E\setminus(-E)\) is
 semi-bounded.\\[1.0ex]
  \hspace*{1.0ex}
Then the set \(E\) is symmetric.
\end{theorem}

\begin{lemma}
\label{dens}
 Let
\(E,\,\,E\subset\mathbb{R}\) be a set of positive measure:
\(\textup{mes}\,(E)>0,\) and the set \(S,\,\,S\subset\mathbb{R}\),
is semi-bounded. Then the set of all functions of the form
\(y(t)=\frac{1}{\sqrt{2\pi}}\int_{E}e^{it\xi}\,x(\xi)\,d\xi,\,
t\in{}S\), where \(x\) runs over \(L^2(E)\), is dense in
\(L^2(S)\).\\
\end{lemma}
\begin{proof}
Assume for definiteness that the set \(S\) is bounded from above,
say \(S\subseteq(-\infty,b]\), where \(b<\infty\). If the set of
all such \(y\) is not dense in \(L^2(S)\), then there exists
\(v\in{}L^2(S),\,v\not=0,\) such that
\(\int\limits_{S}v(t)\overline{y(t)}\,dt=0\) for all \(y(t)\).  From this follow that %
\(\int\limits_{S}v(t)e^{-it\xi}\,dt=0\,\,\forall\,\xi\in{}E\).
Since \(S\subseteq(-\infty,b]\), the function
\(f(\xi)=e^{ib\xi}\int\limits_{S}v(t)e^{-it\xi}\,dt,\,\,\xi\in\mathbb{R}\),
belongs to the Hardy class \(H^2_{+}\). Since \(v\in{}L^2(S)\) is
non-zero, \(f\) is a non-zero function from \(H^2_{+}\). Moreover,
\(f(\xi)=0\) for \(\xi\in{}E\). However, the non-zero function
from the Hardy class can not vanish on the set of positive
measure.
\end{proof}

\begin{remark}
\end{remark}

\begin{proof}[Proof of Theorem \ref{NiS}]
We show that if the set \(S\) is not symmetric, that is if
\(\textup{mes}\,(E\setminus(-E))>0\), then the condition
\eqref{NCoE} is violated for some
\(y(t)=\frac{1}{\sqrt{2\pi}}\int_{E}e^{it\xi}\,x(\xi)\,d\xi\),
where \(x\in{}L^2(E)\). Then by Lemma \ref{NCo}, the operator
\(\mathscr{F}_E\) is not normal.

We first present the proof assuming that the set
\(E\setminus(-E)\) is bounded. If the set \(E\setminus(-E)\) is
bounded, then the set
 \begin{equation}%
 \label{Sym}%
 S\stackrel{\textup{def}}{=}(E\setminus(-E))\cup((-E)\setminus{}E)
 \end{equation}
 is bounded as well. We define the function \(g(t)\) as
 \[g(t)=1,\,t\in(E\setminus(-E)),\,\,g(t)=0,\,t\in((-E)\setminus{}E)\,.\]
 Since the sets \((E\setminus(-E))\) and \(((-E)\setminus{}E)\) do
 not intersect, the definition of the function \(g\) is not contradictory.
 According to \textup{Lemma, \ref{dens}}, for every number \(\varepsilon>0\)
 there exists a function \(y(t)=\frac{1}{\sqrt{2\pi}}\int_{E}e^{it\xi}\,x(\xi)\,d\xi\) such that %
\[\int\limits_{S}\vert{}g(t)-y(t)\vert\,dt\leq\varepsilon^2\,\textup{mes}\,(E\setminus(-E)).\]
Since \(g(t)=1\) on \((E\setminus(-E))\),
\[\int\limits_{E\setminus(-E)}\vert{}y(t)\vert^2\,dt\geq(1-\varepsilon)^2\textup{mes}\,(E\setminus(-E))\,.\]
Since \(g(t)=0\) on \((-E)\setminus{}E)\),
\[\int\limits_{(-E)\setminus{}E}\vert{}y(t)\vert^2\,dt\leq\varepsilon^2\textup{mes}\,(E\setminus(-E))\,.\]
Choosing \(\varepsilon<1/2\), we find
\(y(t)=\frac{1}{\sqrt{2\pi}}\int_{E}e^{it\xi}\,x(\xi)\,d\xi\) for
which the equality \eqref{NCoE} is violated.

If the set \(E\setminus(-E)\) is semi-bounded but not bounded,
then the set \(S\), \eqref{Sym}, is not bounded "in both
directions". We assume for definiteness that set
\(E\setminus(-E)\) is bounded from above, say\footnote{%
Strictly speaking, this condition  should be formulated as %
\(\textup{mes}\,\big((E\setminus(-E))\cap(b,\infty)\big)=0\).
} %
\((E\setminus(-E))\subset(-\infty,b]\), where \(b<\infty\). We
construct such a function \(y(t)\) of the form \eqref{p} for which
\begin{equation}
\label{NCoEi}%
 \int\limits_{E\setminus(-E)}|y(t)|^2dt<\int\limits_{(-E)\setminus{}E}|y(t)|^2dt\,,
\end{equation}
 Since
the set \(E\setminus(-E)\) is bounded from above but not bounded, %
\[\textup{mes}\,\big(((-E)\setminus(E))\cap(b,\infty)\big)>0\,.\]
Therefore there exists a finite interval
\([p,q],\,\,[p,q]\in(b,\infty)\), such that
\[\textup{mes}\,\big({[p,q]\cap((-E)\setminus{}E)}\big)>0\,.\]
Let
\[S=(-\infty,q]\cap\big((E\cap(-E))\cup((-E)\cap{}E)\big),\]
\[\quad g(t)=
\begin{cases}
0,&\text{ if }\quad{}t\in{}S,\ \ -\infty<t<p,\\
1, &\text{ if }\quad{}t\in{}S,\ \ \phantom{--}p\leq{}t\leq{}q\,.
\end{cases}
\]
Clearly,
\begin{equation}%
\label{Ge}%
(E\setminus(-E))\subset{}S\cap(-\infty,p),\,\,S\cap[p,q]=((-E)\setminus{}E)\cap[p,q]\,.
\end{equation}%
and
\begin{equation}
\label{CrIn}%
\int\limits_{E\setminus(-E)}|g(t)|^2dt<\int\limits_{[p,q]\cap((-E)\setminus{}E)}|g(t)|^2dt\,.
\end{equation}
(The left hand side of this inequality is equal to zero, and the
right hand side\,---\,to the strictly positive number
\(\textup{mes}\,\big({[p,q]\cap((-E)\setminus{}E)}).\))\\
 By Lemma \ref{dens}, for any \(\varepsilon>0\) there exists a
function \(y\) of the form \eqref{p} such that
\begin{equation}%
\label{App}%
 \int\limits_S|y(t)-g(t)|^2dt<\varepsilon^2\,.
\end{equation}%
 If
\(\varepsilon\) is small enough, then from
\eqref{Ge}-\eqref{CrIn}-\eqref{App} it follows that
\[\int\limits_{E\setminus(-E)}|y(t)|^2dt<\!\!\!\!\int\limits_{((-E)\setminus{}E)\cap[p,q]}|y(t)|^2dt\,,\]
and all the more the inequality \eqref{NCoEi} holds.
\end{proof}

\begin{theorem}
\label{HSO}%
 Assume that the Lebesgue measure of the set \(E\) is
finite. Then \hspace*{2.0ex}\textup{\textsc{1}}. The operator
\(\mathscr{F}_E^{\,\ast}\mathscr{F}_E\) is an integral operator:
\begin{equation}
\label{11}%
(\mathscr{F}_E^{\,\ast}\mathscr{F}_Ex)(t)=\int\limits_E
K_E(t,s)x(s)\,ds,
\end{equation}
with the kernel
\begin{equation}
\label{12}%
 K_E(t,s)=\int\limits_{E}e^{i\xi(t-s)}\,d\xi,\quad (t,\,s\in{}E)\,.
\end{equation}
\hspace*{2.0ex}\textup{\textsc{2}}. The operator
\(\mathscr{F}_E^{\,\ast}\mathscr{F}_E\) is a trace class operator:
\begin{equation}
\label{13}%
\textup{trace}\,\mathscr{F}_E^{\,\ast}\mathscr{F}_E=(\textup{mes}\,E)^2\,.
\end{equation}
The operator \(\mathscr{F}_E\) belongs to to the class
\(\mathfrak{S}_2\) of Hilbert-Schmidt operators:
\begin{equation}
\label{HSF}%
\Vert\mathfrak{F}_E\Vert_{\mathfrak{S}_2}=\textup{mes}\,E\,.
\end{equation}
 In particular, the operator \(\mathscr{F}_E\) is a compact
operator.\\
 \hspace*{2.0ex}\textup{\textsc{3}}.
 The trace norm of the operator \(\mathscr{F}_E\) satisfy the
 condition
\begin{equation}
\label{HST}%
(\textup{mes}\,E)^2\leq\Vert\mathfrak{F}_E\Vert_{\mathfrak{S}_1}\leq\infty\,.
\end{equation}
\end{theorem}
\proof \ 1.\,The representation \eqref{11}-\eqref{12} is a direct
consequence of the equalities \eqref{6} and \eqref{8} and of the
rule for calculation the kernel of the product of two integral
operators  in terms of their kernels.\\ %
 2.\, The kernel \(
K_E(t,s)\), \eqref{12}, is positive definite, bounded and
uniformly continuous for \((t,s)\in{}E\times{}E.\) As it claimed
in \cite[Chap.3, sect. 10]{GoKr}, from these properties of the
kernel of an integral operator it follows that this operator is a
trace class operator, and that its trace is equal to the integral
\(\int\limits_EK_E(t,t)\,dt\,.\) (See the last paragraph of
section 10 of the quoted reference.)\\
3.\, The equality \eqref{13} means that
\begin{equation}
\label{HSd}%
 \sum\limits_{1\leq_j<\infty}(s_j(\mathscr{F}_E))^2=(\textup{mes}\,E)^2\,,
\end{equation}
where \(s_j(\mathscr{F}_E)\) are the singular value of the
operator \(\mathscr{F}_E\). In view of \eqref{Co},
\[s_j(\mathscr{F}_E)\leq{}1,\quad 1\leq{}j<\infty\,.\]
Thus,
\begin{equation}
\label{Td}
 \sum\limits_{1\leq{}j<\infty}s_j(\mathscr{F}_E)\geq(\textup{mes}\,E)^2\,.
\end{equation}

To study under which conditions the operator \(\mathscr{F}_E\), or
what is the same (see Remark \ref{OtDe}) the operator
\(P_E\mathscr{F}P_E\), belongs
 to the trace class \(\mathfrak{S}_1\),
we have to consider the more general operator
\begin{equation}
\label{CFO}%
 \mathscr{F}_{S_1,S_2}=P_{S_2}\mathscr{F}P_{S_1},
\end{equation}
 where
\(S_1,\,S_2\subset\mathbb{R}\) are measurable sets, and for the
set \(S,\,S\subset\mathbb{R}\), the operator
\(P_S:L_2(\mathbb{R})\to{}L_2(\mathbb{R})\), is defined as
\begin{equation}%
\label{DeP}%
(P_Sx)t)=\chi_S(t)x(t),\,\textup{where}%
\,\,\chi_S(t)=1,\,t\in{}S,\,\,\chi_S(t)=0,\,\,t\notin{}S\,.
\end{equation}%

\begin{theorem}
\label{HSFM}%
 If the truncated Fourier operator \(\mathscr{F}_{E}\)
is a Hilbert-Schmidt operator:
\(\mathscr{F}_{E}\in\mathfrak{S}_2\), then the set \(E\) is of
finite measure, and the equality \eqref{HSF} holds.
\end{theorem}
\begin{proof} The equality \eqref{HSF} was obtained under the
assumption that \(\textup{mes}\,{E}<\infty\). If
\(\textup{mes}\,{E}=\infty\), the kernel \(K_E(t,s)\), \eqref{12},
is not well defined, and the reasoning used in the proof of
Theorem \ref{HSO} is not applicable. Consider the set
\(E_n=E\cap{}[-n,n]\) and the operator \(\mathscr{F}_{E_n}\). The
operator \(\mathscr{F}_{E_n}\) can be identified with the operator
\(P_{E_n}\mathscr{F}_{E}P_{E_n}\) (see Remark \ref{OtDe}\}), where
\(P_{E_n}\) is the orthoprojector in \(L^2(\mathbb{R})\):
\((P_{E_n}x)(t)=\chi_{E_n}(t)x(t)\). Therefore, %
\(\Vert\mathscr{F}_{E_n}\Vert_{\mathfrak{S}_2}\leq\Vert\mathscr{F}_{E}\Vert_{\mathfrak{S}_2}.\)
On the other hand, the set \(E_n\) is of finite measure, and the
formula \eqref{HSF} is applicable to \(E_n\):
\(\textup{mes}\,E_n=\Vert\mathscr{F}_{E_n}\Vert_{\mathfrak{S}_2}\,.\)
Thus, for every \(n\),
\(\textup{mes}\,E_n\leq\Vert\mathscr{F}_{E}\Vert_{\mathfrak{S}_2}\,.\)
Turning \(n\) to infinity, we obtain that
\(\textup{mes}\,E\leq\Vert\mathscr{F}_{E}\Vert_{\mathfrak{S}_2}<\infty.\)
\end{proof}
\begin{lemma}
\label{TEL}%
 Assume that \(S_1,\,S_2\) are bounded measurable sets:
\begin{equation}%
\label{CoB}%
 S_1,S_2\subseteq{}[-R,R],\ \ \ \textup{for some} \ \
R\in(0,\infty).
\end{equation}%

 Then the operator
\(\mathscr{F}_{S_1,S_2}\) belongs to the trace class
\(\mathfrak{S}_1\), and its trace norm
\(\Vert\mathscr{F}_{S_1,S_2}\Vert_{\mathfrak{S}_1}\) admits the
estimate
\begin{equation}
\label{ETR} \Vert\mathscr{F}_{S_1,S_2}\Vert_{\mathfrak{S}_1}\leq
(\textup{mes}\,S_1)^{1/2}\cdot(\textup{mes}\,S_2)^{1/2}\cdot{}e^{R^2}\,.
\end{equation}
\end{lemma}
\begin{proof} The operator \(\mathscr{F}_{S_1,S_2}\) is an integral operator in
the space \(L^2(\mathbb{R})\) with the kernel
\(k(t,\xi)=\chi_{S_2}(t)e^{it\xi}\chi_{S_1}(\xi)\) , which is the
sum of a rank-one kernels: %
\[k(t,\xi)=\sum\limits_{0\leq{}j<\infty}i^jk_j(t\,,\xi),\quad k_j(t\,,\xi)=
\frac{1}{j!}\cdot\chi_{S_2}(t)t^{j}\cdot\xi^{j}\chi_{S_1}(\xi)\,.\] %
The one-dimensional integral operator \(K_j\geq{}0\) with the
kernel \(k_j(t,\xi)\) admits the estimate
\[\Vert{}K_j\Vert_{{\frak S}_1}%
\leq\frac{1}{j!}\cdot\Vert\chi_{S_2}(t)t^{j}\Vert_{L^2(\mathbb{R})}%
\cdot{}\Vert\xi^{j}\chi_{S_1}(\xi)\Vert_{L^2(\mathbb{R})} \,.\]
Since \(S_1\subset[-R,R],\,S_2\subset[-R,R]\),
\[\Vert\chi_{S}(t)t^{j}\Vert_{L^2(\mathbb{R})}\leq(\textup{mes}\,S)^{1/2}R^j,\,\,\,S=S_1,\,S_2\,.\]

 Therefore,
 \[\Vert\mathscr{F}_{S_1,S_2}\Vert_{\mathfrak{S}_1}\leq%
\sum\limits_{0\leq{}j<\infty}\Vert{}K_j\Vert_{\mathfrak{S}_1}\leq
(\textup{mes}\,S_1)^{1/2}\cdot(\textup{mes}\,S_2)^{1/2}\cdot\!\!
\sum\limits_{0\leq{}j<\infty} \frac{R^{2j}}{j!}\,.
 \]
\end{proof}
The estimate \eqref{ETR} shows that if the set \(E\) is bounded,
then the operator \(\mathscr{F}_E\) is a trace class operator.
However this estimate does not work if the set \(E\) is unbounded.

\begin{theorem}%
\label{CTO}%
Let the set \(E,\,E\subset\mathbb{R}\), satisfy the condition
\begin{equation}
\label{CST}%
\sum\limits_{j\in\mathbb{Z}}\,(\textup{mes}\,(E_j))^{1/2}<\infty,\
\end{equation}
where
\begin{equation}
\label{cts}%
 E_j=E\cap[j-1/2,j+1/2]\,,\quad j\in\mathbb{Z}.
\end{equation}
 Then the operator \(\mathscr{F}_E\) is a trace class
operator.
\end{theorem}%
The following lemma is a modification of Lemma \ref{TEL}:
\begin{lemma}
\label{ctc}%
 Let \(E,\,E\subset{}\mathbb{R}\), be a measurable set, and the
 operator
\( \mathscr{F}_{E_p,E_q}=P_{E_q}\mathscr{F}P_{E_p}\) is defined by
\eqref{CFO}-\eqref{DeP}, with \(S_1=E_p,\,S_2=E_q\).

Then
\begin{equation}
\label{BE}%
\Vert\mathscr{F}_{E_p,E_q}\Vert_{\mathfrak{S}_1}\leq%
e^{1/4}\cdot(\textup{mes}\,E_p)^{1/2}\cdot(\textup{mes}\,E_q)^{1/2},\quad
\forall\,p,q\in\mathbb{Z}.
\end{equation}
\end{lemma}
\begin{proof}
For \(p=0,\,q=0\), the estimate \eqref{BE} is the special case of
Lemma \ref{TEL} corresponding \(S_1=E_0,\,S_2=E_0, R=1/2\). The
general case of arbitrary integers \(p\) and \(q\) can be reduced
to the case \(p=0,\,q=0\) by means of translation. The sets
\(S_1=-p+E_p\) and \(S_2=-q+E_q\) are contained in the interval
\([-1/2,1/2]\), and the operator \(\mathscr{F}_{E_p,E_q}\) is
related to the operator \(\mathscr{F}_{S_1,S_2}\) by the equality
\(\mathscr{F}_{E_p,E_q}=U_q\,\mathscr{F}_{S_1,S_2}U_p\), where
\(U_r,\,r\in\mathbb{Z}\), is the unitary operators:
\((U_rx)(t)=e^{irt}x(t)\).
\end{proof}
\begin{proof}[Proof of Theorem \ref{CTO}]
Identifying the operator \(\mathscr{F}_E\) with the operator
\(P_E\mathscr{F}P_E\), we represent it as  the double sum
\[\mathscr{F}_E=\sum\limits_{p\in\mathbb{Z},q\in\mathbb{Z}}P_{E_q}\mathscr{F}P_{E_p}\,,\]
hence
\[\Vert\mathscr{F}_E\Vert_{\mathfrak{S}_1}%
\leq\sum\limits_{p\in\mathbb{Z},q\in\mathbb{Z}}\Vert{}P_{E_q}%
\mathscr{F}P_{E_p}\Vert_{\mathfrak{S}_1}\,.\] Applying to the
summand in the right hand side the estimate \eqref{BE}, we obtain
that
\begin{equation}
\label{FEs}%
 \Vert\mathscr{F}_E\Vert_{\mathfrak{S}_1}\leq{}e^{1/4}\big(\sum\limits_{j\in\mathbb{Z}}%
 (\textup{mes}\,E_j)^{1/2}\big)^{2}\,,
\end{equation}
where \(E_j\) is defined in \eqref{cts}.
\end{proof}
\begin{remark}
Since \(\textup{mes}\,E_j\leq1\), then
\[\textup{mes}\,E=\sum\limits_{j\in\mathbb{Z}}\textup{mes}\,E_j\leq%
\big(\sum\limits_{j\in\mathbb{Z}}%
 (\textup{mes}\,E_j)^{1/2}\big)^{2}\,.\]
 Thus the set \(E\) for which the expression in the right hand side of \eqref{FEs} is
 finite automatically satisfy the condition
 \(\textup{mes}\,E<\infty\). However there are sets \(E\) of
 finite measure for which the expression in the right hand side of
\eqref{FEs} is infinite. For example,
\(E=\bigcup\limits_{1\leq{}j<\infty}[j-j^{-2},j+j^{-2}]\).
\end{remark}

It should be mention that Theorem \ref{CTO} is related to some
results by M.S.Birman and M.Z.Solomyak, \cite[Theorem 11.1]{BiSo},
and may be considered as a special case of their result. However
our presentation is more direct and simple.

Theorem \ref{CTO} is precise:

\begin{theorem}[\textnormal{B.Simon, \cite[Proposition 4.7]{Sim}}]
Let the operator \(\mathscr{F}_E\) be a trace class operator. Then
the set \(E\) satisfy the condition \eqref{CST}\,-\,\eqref{cts}.
\end{theorem}
\begin{proof}
We identify the operator \(\mathscr{F}_E\) with the operator
\(P_E\mathscr{F}P_E\). (See \eqref{CaPr}, \eqref{CaF}, and Remark
\ref{OtDe}.) Since the Fourier-Plancherel operator \(\mathscr{F}\)
is unitary, the operator \(\mathscr{F}^{-1}\mathscr{F}_E=
\mathscr{F}^{-1}P_E\mathscr{F}P_E\) is a trace class operator as
well:
\begin{equation}%
\label{TCO}%
\mathscr{F}^{-1}P_E\mathscr{F}P_E\in\mathfrak{S}_1\,.
\end{equation}%
We are to deduce from \eqref{TCO}, that the set \(E\) satisfy the
condition \eqref{CST}\,-\,\eqref{cts}. According to Theorem
\ref{HSFM}, the condition \eqref{TCO} implies that the set \(E\)
is of finite measure: \(\textup{mes}\,E<\infty\). (See
\eqref{HST}.) Therefore, the function
\begin{equation}%
\label{CoKe}%
h_E(t)=\frac{1}{2\pi}\int\limits_{E}e^{-i\xi{}t}\,d\xi\,,\quad
t\in\mathbb{R}\,,
\end{equation}%
is well defined and continuous on \(\mathbb{R}\). The operator
\begin{equation}%
\label{CoOp}%
 \mathscr{C}_e=\mathscr{F}^{-1}P_E\mathscr{F}P_E
\end{equation}%
 can be
represented as the product of the multiplication and the
convolution operators:
\begin{equation}
(\mathscr{C}_Ex)(t)=\int\limits_{\mathbb{R}}h_E(t-\xi)\,\chi_E(\xi)\,x(\xi)\,d\xi\,.
\end{equation}
We assume that the operator \(\mathscr{C}_E\) is a trace class
operator, acting in the space \(L^2(\mathbb{R})\). We have to
derive from here that the set \(E\) satisfy the condition
\eqref{CST}\,-\,\eqref{cts}.

 The value \(h(0)=\frac{1}{2\pi}\textup{mes}\,E\) is strictly positive, and the function
\(h\) is continuous.
  Therefore there
exists \(\delta>0\) such that
\begin{equation}
\label{ERP}%
\textup{Re}\,h(t)>\frac{1}{4\pi}\,\textup{mes}\,E,\quad
t\in[-\delta,\,\delta]\,.
\end{equation}
 Since the operator \(\mathscr{C}_E\) is a trace class operator
in \(L^2(\mathbb{R})\), for any two orthonormal systems
\(\lbrace\varphi_m\rbrace_{m\in{}M}\) and
\(\lbrace\psi_m\rbrace_{m\in{}M}\), \(M\subseteq\mathbb{Z}\) is an
indexing set, the inequality holds
\begin{equation}
\label{TCE}
\sum\limits_{m\in{}M}|\langle\mathscr{C}_E\varphi_m,\psi_m%
\rangle_{L^2_E}|\leq\Vert\mathscr{C}_E\Vert_{\mathfrak{S}_1}<\infty\,.
\end{equation}
We obtain the information concerning the set \(E\) choosing the
systems \(\lbrace\varphi_j\rbrace_{m\in{}M}\) and
\(\lbrace\psi_m\rbrace_{j\in{}M}\) by an appropriate way.

 Let
\begin{equation}
\label{SbD}
Q_{m,\delta}=[m\delta-\delta/2,m\delta+\delta/2)\,,\quad
m\in\mathbb{Z}\,.
\end{equation}
The intervals \(Q_{m,\delta},\,m\in\mathbb{Z}\), form the
partition of the real axis. Let
\begin{equation}%
\label{PaE}%
 E_{m,\delta}=E\cap{}Q_{m,\delta}\,,\quad m\in\mathbb{Z}\,.
\end{equation}%
We will prove that for chosen \(\delta\),
\begin{equation}%
\label{L12d}%
\sum\limits_{m\in\mathbb{Z}}(\textup{mes}\,E_{m,\delta})^{1/2}<\infty\,.
\end{equation}%
Let
\(M=\lbrace{}m\in\mathbb{Z}:\,\textup{mes}\,E_{m,\delta}>0\rbrace.\)
For \(m\in{}M\), we set
\begin{subequations}%
\label{CSB}%
\begin{align}
\label{CSBa}%
\varphi_m(t)&=(\textup{mes}\,Q_{m,\delta})^{-1/2}\cdot%
\chi_{_{\scriptstyle{}Q_{m,\delta}}}(t),\\ %
\label{CSBb}%
\psi_m(t)&=(\textup{mes}\,E_{m,\delta})^{-1/2}\cdot%
\chi_{_{\scriptstyle{}E_{m,\delta}}}(t)\,.
\end{align}
\end{subequations}%
The systems \(\lbrace\varphi_m\rbrace_{m\in{}M}\) and
\(\lbrace\psi_m\rbrace_{m\in{}M}\), defined by \eqref{CSB}, are
orthonormal. Let us calculate end estimate the scalar product %
\(\langle\mathscr{C}_E\,\varphi_m,\psi_m\rangle_{{}_{\!L^2(\mathbb{R})}}\).
According to \eqref{CoOp} and \eqref{CSB},
\begin{multline}%
\label{ScP}%
\langle\mathscr{C}_E\,\varphi_m,\psi_m\rangle_{{}_{\!L^2(\mathbb{R})}}=\\
(\textup{mes}\,Q_{m,\delta})^{-1/2}(\textup{mes}\,E_{m,\delta})^{-1/2}
\iint\limits_{\substack{t\in{}Q_{m,\delta}\\
\xi\in{}Q_{m,\delta}}}h(t-\xi)\chi_{_{\scriptstyle{}E_{m,\delta}}}(\xi)\,d\xi\,.%
\end{multline}
According to \eqref{SbD}, for
\(t\in{}Q_{m,\delta},\,\xi\in{}Q_{m,\delta}\), the inequality
\(\vert{}t-\xi\vert\leq{}\delta\) holds. Together with
\eqref{ERP}, this yields
\begin{equation}%
\label{EhB}%
 \textup{Re}\,h(t-\xi)\geq c>0\ \ \textup{for}\ \ %
 t\in{}Q_{m,\delta},\,\xi\in{}Q_{m,\delta},\ \textup{where} \ %
 c=\frac{1}{4\pi}\textup{mes}\,E\,.
\end{equation}%
Invoking \eqref{ScP}, we obtain
\begin{multline}%
\label{BES}%
\vert\langle\mathscr{C}_E\,\varphi_m,\psi_m\rangle_{{}_{\!L^2(\mathbb{R})}}\vert\geq
\textup{Re}\,\langle\mathscr{C}_E\,\varphi_m,\psi_m\rangle_{{}_{\!L^2(\mathbb{R})}}=\\
(\textup{mes}\,Q_{m,\delta})^{-1/2}(\textup{mes}\,E_{m,\delta})^{-1/2}
\iint\limits_{\substack{t\in{}Q_{m,\delta}\\
\xi\in{}Q_{m,\delta}}}\textup{Re}\,h(t-\xi)\chi_{_{\scriptstyle{}E_{m,\delta}}}(\xi)\,d\xi\,.
\end{multline}%
Finally, taking into account that
\(\textup{mes}\,Q_{m,\delta}=\delta\), we get
\begin{equation}%
\label{BESf}%
\vert\langle\mathscr{C}_E\,\varphi_m,\psi_m\rangle_{{}_{\!L^2(\mathbb{R})}}\vert\geq
c\,\delta\,(\textup{mes}\,E_{m,\delta})^{1/2}\,.
\end{equation}%
From here and \eqref{TCE}, the condition \eqref{L12d} follows.

The condition \eqref{L12d} is almost what we need. We obtain the
the condition for any \(\delta\) satisfying the condition
\eqref{ERP}. We need the condition \eqref{L12d} for \(\delta=1\).
(For \(\delta=1\), this is the condition \eqref{CST}.)

Actually, if the the condition \eqref{L12d} is fulfilled for
\emph{some} positive \(\delta\), then it is fulfilled for
\emph{any} positive \(\delta\). We show this in the generality
which we need. Because we may diminish \(\delta\) without to
violate the condition \eqref{ERP}, we choose \(\delta\) of the
form
\begin{equation}
\label{CCD}%
\delta=\frac{1}{N}\,, \quad N \ \textup{is a positive integer
which is large enough\,.}
\end{equation}
If \(a_1,\,a_2,\,\ldots\,,\,a_N\) are non-negative numbers, then
\[(a_1+a_2+\,\cdots\,+a_N)^{1/2}\leq
\sum\limits_{1\leq{}k\leq{}N}a_k^{1/2}\,.\] It is clear that for
\(\delta=1/N\), either the sets \(Q_{j,1}\) and \(Q_{m,\delta}\)
 do not intersect, or the set \(Q_{m,\delta}\) is contained
in \(Q_{j,1}\). Moreover, the total number of the sets
\(Q_{m,\delta}\) which are contained in \(Q_{j,1}\) is equal to
\(N\). Thus,
\[\textup{mes}\,(E\cap{}Q_{j,1})=\sum\limits_{m:\,Q_{m,\delta}\subset{}Q_{j,1}}\textup{mes}\,(E\cap{}Q_{m,\delta})\,,\]
and the sum in the right hand side contains precisely \(N\)
summands. Therefore, for every \(j\in\mathbb{Z}\),
\[(\textup{mes}\,(E\cap{}Q_{j,1}))^{1/2}\leq{}
\sum\limits_{m:\,Q_{m,\delta}\subset{}Q_{j,1}}(\textup{mes}\,(E\cap{}Q_{m,\delta}))^{1/2}\,,\]
and
\begin{equation}%
\label{EN}
\sum\limits_{j\in\mathbb{Z}}(\textup{mes}\,(E\cap{}Q_{j,1}))^{1/2}\leq{}
\sum\limits_{m\in\mathbb{Z}}(\textup{mes}\,(E\cap{}Q_{j,\delta}))^{1/2}\,.
\end{equation}%
\end{proof}

\begin{question} For which sets \(E\) the operator \(\mathscr{F}_E\)
is compact?
\end{question}

 In the next part of this work we  embbark on a more
detail discussion of spectral properties of the operators
\(\mathscr{F}_E\) in three important  cases:\\
\hspace*{1.0ex}\(\boldsymbol{\cdot}\)  \(E=\mathbb{R}\);\\ %
\hspace*{1.0ex}\(\boldsymbol{\cdot}\)  \(E\) is an arbitrary
symmetric finite
interval: \(E=[-a,a]\), \(a\in]0,\infty[\); \\ %
\hspace*{1.0ex}\(\boldsymbol{\cdot}\)   \(E=[0,\infty]\).

The case \(E=\mathbb{R}\) has already been studied in great
details. We review shortly the main facts about this case. Also
the case \(E=[-a,a],\,\,a\in]0,\infty[\), was  already considered.
However this case is more complicated than previous, and some
questions remain open. To the best of our knowledge, the case
\(E=[0,\infty[\) was not studied until now.


\vspace{4.0ex}
\begin{minipage}[h]{0.45\linewidth}
Victor Katsnelson\\[0.2ex]
Department of Mathematics\\
The Weizmann Institute\\
Rehovot, 76100, Israel\\[0.1ex]
e-mail:\\
{\small\texttt{victor.katsnelson@weizmann.ac.il}}
\end{minipage}

\vspace{3.0ex}
\begin{minipage}[h]{0.45\linewidth}
Ronny Machluf\\[0.2ex]
Department of Mathematics\\
The Weizmann Institute\\
Rehovot, 76100, Israel\\[0.1ex]
e-mail:\\
\texttt{ronny-haim.machluf@weizmann.ac.il}
\end{minipage}
\end{document}